\documentclass[final,onefignum,onetabnum]{siuro250211}

\usepackage{lipsum}
\usepackage{amsfonts}
\usepackage{graphicx}
\usepackage{epstopdf}
\usepackage{algorithmic}
\ifpdf
  \DeclareGraphicsExtensions{.eps,.pdf,.png,.jpg}
\else
  \DeclareGraphicsExtensions{.eps}
\fi

\newcommand{\inmag}[1]{\mbox{\color{magenta} #1}}

\makeatletter
\renewcommand*\env@matrix[1][\arraystretch]{%
  \edef\arraystretch{#1}%
  \hskip -\arraycolsep
  \let\@ifnextchar\new@ifnextchar
  \array{*\c@MaxMatrixCols c}}
\makeatother

\usepackage{enumitem}
\setlist[enumerate]{leftmargin=.5in}
\setlist[itemize]{leftmargin=.5in}

\newsiamremark{remark}{Remark}
\newsiamremark{hypothesis}{Hypothesis}
\crefname{hypothesis}{Hypothesis}{Hypotheses}
\newsiamthm{claim}{Claim}
\newsiamremark{fact}{Fact}
\crefname{fact}{Fact}{Facts}

\def\reals{\mathbb R}

\title{A Modified BGK Collision Operator for Exact Conservation \\ in Numerical Solutions of Boltzmann-BGK }

\author{Vienna B. Rossmanith\thanks{Ames High School, Ames, IA, USA 
  (\email{887432ros@ames.k12.ia.us})}}

\dedication{\small\textit{Project advisor: Preeti Sar \thanks{Fusion Energy Division, Oak Ridge National Laboratory, Oak Ridge, TN, USA}}}

\usepackage{amsopn}

\begin{document}

\maketitle

\begin{abstract}
Ideal gases can be modeled by the Boltzmann equation from statistical physics. Instead of trying to track the position and velocity of a large number of gas molecules, it is possible to describe the particles with a particle distribution function. The Boltzmann equation provides the rule for evolving the distribution function over time, allowing one to simulate the gas dynamics. In this work, we develop a novel numerical method for solving the 1D1V Boltzmann-BGK equation. Several important ingredients are combined to create an accurate, efficient, robust numerical method valid in all flow regimes. First, we make use of operator splitting to create separate transport and collision sub-steps, each of which is easier to discretize than the whole system; second, we introduce a third-order accurate Lax-Wendroff-type scheme for the transport sub-step; third, we make use of the second-order L-stable TR-BDF method for the collision sub-step. The final component we introduce is a novel treatment of the collision sub-step to guarantee that the total mass, momentum, and energy are conserved even in the presence of a truncated velocity range and quadrature errors in computing moments. The key idea is to multiply the Maxwell-Boltzmann distribution in the collision sub-step by a quadratic Hermite polynomial, where the coefficients in the polynomial are chosen to ensure exact conservation. The resulting scheme is verified on piecewise constant initial data with periodic boundary conditions; exact conservation up to machine precision is demonstrated for mass, momentum, and energy. The resulting method is implemented in a freely available Java code with Python plotting routines.
\end{abstract}


\section{Introduction}
A gas is a collection of unbound molecules; in the absence of external forces, each molecule moves with some velocity until it collides with another molecule. The collisions act to drive the gas towards a thermodynamic equilibrium. Mathematically, one could consider a range of models, including (1) molecular dynamics, which is valid at very small scales, where each molecule in some small portion of the gas is tracked; (2) kinetic Boltzmann equations, which is valid at intermediate scales, where a particle distribution function replaces individual molecules; and (3) fluid models such the Navier-Stokes equation, which are valid at large scales, where macroscopic quantities such as mass, momentum, and energy are tracked (e.g., see \cite{Roohi2025}).

The focus of this work is the kinetic Boltzmann equation; in particular, we consider the 1D1V Boltzmann-BGK equation \cite{Bhatnagar1954ASystems}, in which the full Boltzmann collision operator is replaced by a far simpler collision operator that nonetheless captures key dynamics. While there exist several approaches to solving the Boltzmann-BGK equation, we consider in this work so-called {\it direct kinetic} methods, rather than stochastic methods such as Direct Simulation Monte-Carlo (DSMC) (e.g., see \cite{Pareschi2001AnEquation}). In a direct kinetic approach, various numerical discretizations can be applied on a computational mesh that covers a space/velocity domain that is truncated in velocity.

The goal of this work is to develop a finite difference method for the one-dimensional Boltzmann-BGK model of gas dynamics that is similar in spirit to existing approaches (e.g., Pieraccini and Puppo \cite{Pieraccini2007Implicit-explicitEquations}, Li et al. \cite{article:LiFang2022}, and Sar \cite{thesis:Sar2024}), but with the additional feature that mass, momentum, and energy are exactly conserved. Several important ingredients are combined to create an accurate, efficient, robust numerical method valid in all flow regimes. First, we make use of operator splitting to create separate transport and collision sub-steps, each of which is easier to discretize than the whole system; second, we introduce a third-order accurate Lax-Wendroff-type scheme for the transport sub-step; third, we make use of the second-order L-stable TR-BDF method for the collision sub-step. The final component we introduce is a novel treatment of the collision sub-step to guarantee that the total mass, momentum, and energy are conserved even in a truncated velocity range and quadrature errors in computing moments. The key idea is to multiply the Maxwell-Boltzmann distribution in the collision sub-step by a quadratic Hermite polynomial, where the coefficients in the polynomial are chosen to ensure exact conservation. 

The remainder of this paper is organized as follows. In \S\ref{sec:equations} we briefly describe the 1D1V Boltzmann-BGK equation. In \S\ref{sec:numerical-methods} we develop the main components of the proposed numerical method using operator splitting and finite differences. The scheme described in \S\ref{sec:numerical-methods} does not exactly conserve mass, momentum, and energy. Therefore, in \S\ref{sec:modified-collision} we develop a modified collision operator that is able to recover exact conservation; in this same section, we present the final proposed numerical method in two short algorithms. A numerical example illustrating the benefits of the proposed modification is provided in \ref{sec:numerical-results}. We conclude in \S{sec:conclusion}.


\section{Boltzmann-BGK}
\label{sec:equations}
A particle distribution function (PDF) of time, position, and velocity can be used to describe the movement of molecules in a gas. In a one-dimensional system, the PDF is a function of the following form:
\begin{equation}
\label{eqn:PDF}
    f(t,x,v) : \reals_{\ge 0} \times \reals \times \reals \mapsto \reals_{\ge 0},
\end{equation}
where $t\in \reals_{\ge 0}$ is time, $x\in \reals$ is the space coordinate, and $v\in \reals$ is the velocity coordinate. The total number of particles that can be found in some spatial region $x \in \left[x_1,x_2\right]$ with some range of velocities $v\in \left[v_1, v_2 \right]$ can be computed by integrating the PDF over these values:
\begin{equation}
\text{total number of particles in} \, \left[x_1,x_2\right] \times \left[v_1,v_2\right] = 
\int_{x_1}^{x_2} \int_{v_1}^{v_2} f\left(t,x,v\right) \, dv \, dx.
\end{equation}

According to gas kinetic theory, the PDF satisfies the Boltzmann equation, which in one dimension can be written as follows:
\begin{equation}
\label{eqn:Boltzmann-generic}
    \frac{\partial f}{\partial t} + v \, \frac{\partial f}{\partial x} = \frac{1}{\varepsilon} C(f),
\end{equation}
which consists of two key pieces: (1) the motion of particles (i.e., gas molecules) according to their velocity, $v$, and the collision of particles as modeled by the collision operator, $C(f)$. The parameter $\varepsilon>0$ is called the Knudsen number and is defined as the ratio of the mean-free path (i.e., how far a particle travels on average before undergoing a collision) to an important characteristic length scale in the problem (i.e., size of the domain). Small $\varepsilon>0$ means high collisionality, while large $\varepsilon>0$ means low collisionality.

The net effect of the collision operator on the dynamics is that it drives the gas towards a {\it thermodynamic equilibrium}, which in one dimension is given by the following distribution:
\begin{equation}
\label{eqn:Maxwellian}
    {\mathcal M}\left(t,x,v\right) = \frac{\rho(t,x)}{\sqrt{2\pi \, T(t,x)}} \exp\left[-\frac{\left(v-u(t,x)\right)^2}{2T(t,x)}\right];
\end{equation}
this distribution is referred to as the {\it Maxwell-Boltzmann} distribution and depends on three macroscopic (or fluid) variables:
\begin{align}
\label{eqn:density}
   \text{Mass density:} & \quad \rho(t,x) = \int_{-\infty}^{\infty} f(t,x,v) \, dv = \int_{-\infty}^{\infty} {\mathcal M}(t,x,v) \, dv, \\
\label{eqn:velocity}
    \text{Fluid velocity:} & \quad u(t,x) = \frac{1}{\rho(t,x)}\int_{-\infty}^{\infty} v f(t,x,v) \, dv = 
    \frac{1}{\rho(t,x)}\int_{-\infty}^{\infty} v {\mathcal M}(t,x,v) \, dv, \\
    \begin{split}
\label{eqn:temperature}
        \text{Temperature:} & \quad  T(t,x) = \frac{1}{\rho(t,x)}\int_{-\infty}^{\infty}
    {\bigl(v-u(t,x)\bigr)^2} f(t,x,v) \, dv \\
    & \hspace{15mm} = \frac{1}{\rho(t,x)}\int_{-\infty}^{\infty} {\bigl(v-u(t,x)\bigr)^2} {\mathcal M}(t,x,v) \, dv.
    \end{split}
\end{align}
These variables are also related to two other important quantities:
\begin{align}
\label{eqn:momentum}
    \text{Momentum:} & \quad m(t,x) = \int_{-\infty}^{\infty} v f(t,x,v) \, dv = 
    \int_{-\infty}^{\infty} v {\mathcal M}(t,x,v) \, dv, \\
\label{eqn:energy}
        \text{Energy:} & \quad  {\mathcal E}(t,x) = \int_{-\infty}^{\infty}
    v^2 \, f(t,x,v) \, dv =\int_{-\infty}^{\infty} v^2 \, {\mathcal M}(t,x,v) \, dv,
\end{align}
where
\begin{equation}
m = \rho u \qquad \text{and} \qquad
{\mathcal E} = \rho \left( u^2 + T \right).
\end{equation}

One difficulty in solving equation \eqref{eqn:Boltzmann-generic} is that the collision operator, $C(f)$, can be quite complicated and involve computational expensive integration over phase space. To reduce this complexity, a standard approach is to replace the full Boltzmann operator by something simpler; a commonly used model is the Bhatnagar, Gross, and Krook (BGK) \cite{Bhatnagar1954ASystems} approximation, which mimics the net effect of the Boltzmann collision operator, namely driving the distribution function, $f$, towards its thermodynamic equilibrium: \eqref{eqn:Maxwellian}. The Boltzmann-BGK equation can be written as follows:
\begin{equation}
\label{eqn:Boltzmann-BGK}
    {f}_{,t} + v \, {f}_{,x} = \frac{1}{\varepsilon} C(f) = \frac{1}{\varepsilon} \left( {\mathcal M} - f \right).
\end{equation}
By construction, the mass \eqref{eqn:density}, momentum \eqref{eqn:momentum}, and energy \eqref{eqn:energy} are invariants of the BGK collision operator:
\begin{align}
\label{eqn:cons-mass}
\text{Mass:} & \quad 
\frac{1}{\varepsilon} \int_{-\infty}^{\infty} C(f) \,  dv = \frac{1}{\varepsilon} \int_{-\infty}^{\infty} \left({\mathcal M}-f\right) \, dv =
\frac{1}{\varepsilon} \left( \rho-\rho \right) = 0, \\
\label{eqn:cons-mom}
\text{Momentum:} & \quad 
\frac{1}{\varepsilon} \int_{-\infty}^{\infty} v \, C(f) \,  dv = \frac{1}{\varepsilon} \int_{-\infty}^{\infty} v \left({\mathcal M}-f\right) \, dv =
\frac{1}{\varepsilon} \left( \rho u -\rho u \right) = 0, \\
\label{eqn:cons-energy}
\text{Energy:} & \quad 
\frac{1}{\varepsilon} \int_{-\infty}^{\infty} v^2 \, C(f) \, dv = \frac{1}{\varepsilon} \int_{-\infty}^{\infty} v^2 \left({\mathcal M}-f\right) \, dv =
\frac{1}{\varepsilon} \left( {\mathcal E} - {\mathcal E} \right) = 0.
\end{align}
Note that the fluid velocity \eqref{eqn:velocity} and temperature \eqref{eqn:temperature} are also invariants of the collision operator.

\section{Numerical Method}
\label{sec:numerical-methods}
To approximate the particle distribution function, we consider a domain with finite extent in both $x$ and $v$:
\begin{equation}
\left(x, v\right)  \in  \left[a, b\right] \times \left[v_{\text{min}}, v_{\text{max}} \right].
\end{equation}
The finite extent in $v$ may lead to errors in $f$ (i.e., the Maxwell-Boltzmann distribution \eqref{eqn:Maxwellian} only becomes zero as $|v|\rightarrow \infty$), but is necessary for a mesh-based approach. We then introduce a mesh with uniform rectangles of size $\Delta x$ by $\Delta v$, centered at $(x_i,v_j)$:
\begin{equation}
x_i = a + \left(i-\frac{1}{2}\right) \Delta x, \qquad
v_j = v_{\text{min}} + \left(j-\frac{1}{2}\right) \Delta v,
\end{equation}
for $i = 1,\ldots,N_x$ and $j = 1,\ldots,N_v$, with grid spacings given by
\begin{equation}
\Delta x = \frac{b-a}{N_x}, \qquad \Delta v = \frac{v_{\text{max}}-v_{\text{min}}}{N_v}.
\end{equation}
Note that $N_x \times N_v$ is the total number of grid cells. 

On this mesh, we begin with the following initial condition at $t=0$:
\begin{equation}
	f^{0}_{ij} := f\left(t=0, \, x_i, \, v_j \right) \quad \text{for} \quad i=1,\ldots,N_x \quad \text{and} \quad
		j=1,\ldots,N_v.
\end{equation}
We aim to advance the initial distribution from $t=0$ to some final time $t=T_{\text{final}}$. We will advance the initial $f^{0}$ at $t=0$ to $f^{1}$ at $t=\Delta t$, then advance $f^{1}$ at $t=\Delta t$ to 
$f^{2}$ at $t=2\Delta t$, and so on, until some final value $f^{N}$ at time $t=N \Delta t$,
where $N$ is the total number of time steps with a time step of $\Delta t = T_{\text{final}}/N$.

Recall that the Boltzmann-BGK equation can be written as \eqref{eqn:Boltzmann-BGK}. To simplify the system, we first apply Strang operator splitting \cite{article:St68} to split a single time step of length $\Delta t$ of the numerical method into transport and collision parts:
\begin{equation}
\label{eqn:splitting}
\begin{cases}
\frac{\Delta t}{2}:& \frac{\partial f}{\partial t} = \frac{1}{\varepsilon} \left( {\mathcal M} - f \right), 
\hspace{4mm} \text{\inmag{\tt // Half time-step of collisions only}} \\[6pt]
{\Delta t}:&  \frac{\partial f}{\partial t} + v \, \frac{\partial f}{\partial x} = 0, 
\hspace{8.5mm} \text{\inmag{\tt // Full time-step of transport only}}\\[6pt]
\, \frac{\Delta t}{2}:& \frac{\partial f}{\partial t} = \frac{1}{\varepsilon} \left( {\mathcal M} - f \right). 
\hspace{4.5mm} \text{\inmag{\tt // Half time-step of collisions only}}
\end{cases}
\end{equation}
We describe the numerical discretizations used for both the transport and collision sub-steps below.

\subsection{Transport Sub-Step}
We discretize the transport sub-step,
\begin{equation}
\frac{\partial f}{\partial t} + v \, \frac{\partial f}{\partial x} = 0,
\end{equation}
by a third-order variant of the Lax-Wendroff method \cite{Lax1960SystemsLaws}.
Consider the following Taylor series in time, comparing the PDF at some time $t=t^n$ to the PDF at the time $t=t^n + \Delta t$:
\begin{equation}
\label{eqn:taylor-time}
f\left(t^n + \Delta t, \, x, \, v \right) = \left( f + \Delta t \, 
\frac{\partial f}{\partial t} 
+ \frac{\Delta t^2}{2}  \, \frac{\partial^2 f}{\partial t^2} 
+ \frac{\Delta t^3}{6}  \, \frac{\partial^3 f}{\partial t^3} \right)\Biggl|_{\left(t^n, \, x, \, v \right)} + \, {\mathcal O}\left(\Delta t^4\right).
\end{equation}
Each time derivative can then be replaced by a spatial derivative via the transport equation:
\begin{equation}
\frac{\partial f}{\partial t} = - v \, \frac{\partial f}{\partial x}, \quad
\frac{\partial^2 f}{\partial t^2} = - v \, \frac{\partial}{\partial x} \left( 
     \frac{\partial f}{\partial t} \right) = v^2 \frac{\partial^2 f}{\partial x^2}, \quad
     \frac{\partial^3 f}{\partial t^3} = v^2 \, \frac{\partial^2}{\partial x^2} \left( 
     \frac{\partial f}{\partial t} \right) = -v^3 \, \frac{\partial^3 f}{\partial x^3}.
\end{equation}
Plugging this into \eqref{eqn:taylor-time} and evaluating at the cell center $\left(x_i, v_j \right)$ yields:
\begin{equation}
\label{eqn:taylor-time-2}
f\left(t^n + \Delta t, \, x_i, \, v_j \right) \approx \left( f - v \, \Delta t \, 
\frac{\partial f}{\partial x} 
+ \frac{v^2 \, \Delta t^2}{2}  \, \frac{\partial^2 f}{\partial x^2} 
- \frac{v^3 \, \Delta t^3}{6}  \, \frac{\partial^3 f}{\partial x^3} \right)\Biggl|_{\left(t^n, \, x_i, \, v_j \right)}.
\end{equation}
Finally, all spatial derivatives are replaced by finite differences. For all negative velocity values $\left(v_j<0\right)$, we compute the first, second, and third spatial derivatives at $x_i$ using the highest-order approximation possible using the points $\left( x_{i-1}, \, x_{i}, \, x_{i+1}, \, x_{i+2} \right)$. For all positive velocity values $\left(v_j>0\right)$, we compute the first, second, and third spatial derivatives at $x_i$ using the highest-order approximation possible using the points $\left( x_{i-2}, \, x_{i-1}, \, x_{i}, \, x_{i+1} \right)$. We omit the details, but summarize the resulting update formulas here:
\begin{gather}
\label{eqn:transport1}
  \begin{split}
	\nu_j < 0: \quad f_{ij}  &=  
f^{\star}_{ij} - \frac{\nu_j}{6}  \biggl( 
	-2f^{\star}_{i-1 \, j} - 3 f^{\star}_{ij} + 6 f^{\star}_{i+1 \, j} - f^{\star}_{i+2 \, j}   \biggr) \\
	&+ \frac{\nu_j^2}{2} \biggl( f^{\star}_{i-1 \, j} - 2f^{\star}_{ij} + f^{\star}_{i+1 \, j} \biggr)
	- \frac{\nu_j^3}{6} \biggl( -f^{\star}_{i-1 \, j} + 3f^{\star}_{ij} - 3f^{\star}_{i+1 \, j} + f^{\star}_{i+2 \, j} \biggr),
\end{split} \\
\label{eqn:transport2}
	\begin{split}
	\nu_j > 0: \quad f_{ij}  &=  
f^{\star}_{ij} - \frac{\nu_j}{6}  \biggl( 
	f^{\star}_{i-2 \, j} - 6 f^{\star}_{i-1 \, j} + 3 f^{\star}_{ij} + 2f^{\star}_{i+1 \, j} \biggr) \\
	&+ \frac{\nu_j^2}{2} \biggl( f^{\star}_{i-1 \, j} - 2f^{\star}_{ij} + f^{\star}_{i+1 \, j} \biggr)
	- \frac{\nu_j^3}{6} \biggl( -f^{\star}_{i-2 \, j} + 3f^{\star}_{i-1 \, j} - 3f^{\star}_{ij} + f^{\star}_{i+1 \, j} \biggr),
	\end{split}
\end{gather}
where $\nu_j =  {v_j \, \Delta t}/{\Delta x}$, and in the above formulas, $f^{\star}$ represents the old value of the distribution (i.e., at time $t=t^n$), while $f$ represents the new value (i.e., at time $t=t^n + \Delta t$).

\subsection{Collison Sub-Step}
The collision sub-step requires us to approximately solve
\begin{equation}
\frac{\partial f}{\partial t} = \frac{1}{\varepsilon} \left( {\mathcal M} - f \right).
\end{equation}
We first need to be able to compute the discrete Maxwell-Boltzmann distribution \eqref{eqn:Maxwellian}, which first requires the computation of the mass density \eqref{eqn:density},
fluid velocity \eqref{eqn:velocity}, and temperature \eqref{eqn:temperature}. At the discrete level, we accomplish this via simple Riemann sums:
\begin{gather}
\label{eqn:discrete_moments}
	\begin{pmatrix}[1.2] \rho_i  \\ m_i \\ {\mathcal E}_i \end{pmatrix} = \Delta v \sum_{j=1}^{N_v} \begin{pmatrix}[1.2]
	 1 \\ v_j \\ v_j^2 \end{pmatrix} f_{ij}, \quad
	u_i = \frac{m_i}{\rho_i}, \quad T_i = \frac{{\mathcal E}_i}{\rho_i} -  \left(u_i\right)^2,
\end{gather}
and then evaluate the Maxwell-Boltzmann distribution at each cell center:
\begin{gather}
    \mathcal{M}_{ij} = \frac{\rho_i}{\sqrt{2\pi T_i}} \exp\left[-\frac{\mu_{ij}^2}{2} \right], \quad
    \text{where} \quad \mu_{ij} = \frac{v_j-u_i}{\sqrt{T_i}}.
\end{gather}

To get stable and accurate solutions from the collision sub-step, we make use of the L-stable TR-BDF2 scheme \cite{article:TRBDF2}. The L-stable property allows for the resulting numerical method to be stable for fixed $\Delta t$, $\Delta x$, and $\Delta v$ for any Knudsen number, $\varepsilon>0$, including in the singular limit: $\varepsilon \rightarrow 0^+$. This makes the overall scheme {\it asymptotic-preserving} (AP) \cite{Caflisch1997UniformlyRelaxation,article:Jin1995,Jin1999EfficientEquations}. For a generic ordinary differential equation of the form: $u' = g(u)$, the TR-BDF2 scheme is a 2-stage implicit Runge-Kutta method that can be written as follows:
\begin{equation}
\label{eqn:generic-trbdf2}
U^{\star} = U^n + \frac{\Delta t}{4} \left( g\left(U^n\right) + g\left( U^{\star} \right) \right), \qquad
3 U^{n+1} - 4 U^{\star} + U^n = \Delta t \, g\left(U^{n+1}\right).
\end{equation}

The collision sub-system over a single time-step from $t=t^n$ to $t=t^n + \Delta t$ for the Boltzmann-BGK equation will have the following form:
\begin{equation}
\label{eqn:collision-sub-onestep}
\frac{\partial f}{\partial t} = \frac{1}{\varepsilon} \left( {\mathcal M}^n - f \right),
\end{equation}
where we note that over the time step, $t \in \left[t^n, t^n + \Delta t\right]$, the Maxwell-Boltzmann distribution is constant in time: ${\mathcal M}^n$. Why is it constant in time? Because, as already shown above, the mass \eqref{eqn:cons-mass}, momentum \eqref{eqn:cons-mom}, and energy \eqref{eqn:cons-energy}, are all collision invariants, which also means the Maxwell-Boltzmann distribution is a collision invariant. With a fixed ${\mathcal M}^n$, it is quite simple to apply the TR-BDF2 scheme \eqref{eqn:generic-trbdf2} to \eqref{eqn:collision-sub-onestep}:
\begin{equation}
f^{\star}_{ij} = f^n_{ij} + \frac{\Delta t}{4\varepsilon} \left( 2{\mathcal M}^n_{ij} - f^n_{ij} - f^{\star}_{ij} \right), \qquad
3 f^{n+1}_{ij} - 4 f^{\star}_{ij} + f^n_{ij} = \frac{\Delta t}{\varepsilon} \, \left( {\mathcal M}^n_{ij} - f^{n+1}_{ij}\right).
\end{equation}
After some algebra (i.e., first solve for $f^{\star}$ in terms of $f^{n}$ and ${\mathcal M}^{n}$ in the first stage, then use this result to solve for $f^{n+1}$ in terms of $f^{n}$ and ${\mathcal M}^{n}$ in the second stage), this can be reduced to the following update:
\begin{equation}
\label{eqn:collision-update}
  f^{n+1}_{ij} = \theta {\mathcal M}^n_{ij} + \left( 1 - \theta \right) f^n_{ij}, \quad \text{where} \quad
	\theta = \frac{\Delta t \left(\Delta t + 12 \varepsilon\right)}{\left(\Delta t + 3 \varepsilon\right) \left(\Delta t + 4 \varepsilon\right)}.
\end{equation}
Note that in the current discussion, we are using a time step of size $\Delta t$, while in the splitting paradigm of \eqref{eqn:splitting}, we take half-time steps of the collision sub-step. To take a half-time step, we just need to replace $\Delta t$ in \eqref{eqn:collision-update} with $\Delta t/2$.

\section{Modified Collision Operator}
\label{sec:modified-collision}
As described above in \eqref{eqn:cons-mass}, \eqref{eqn:cons-mom}, and \eqref{eqn:cons-energy}, the mass, momentum, and energy are invariants of the collision operator. At the discrete level, this is only approximately true \cite{Mieussens2000DiscreteDynamics,Pieraccini2007Implicit-explicitEquations} . Consider the following collision sub-step as described above:
\begin{equation}
  f^{n+1}_{ij} = \theta {\mathcal M}^n_{ij} + \left( 1 - \theta \right) f^n_{ij},
\end{equation}
and compute the discrete mass, momentum, and energy via \eqref{eqn:discrete_moments}:
\begin{equation}
  \sum_{j=1}^{N_v} \begin{pmatrix}[1.2] 1  \\ v_j \\ v_j^2 \end{pmatrix} f^{n+1}_{ij} = \theta \sum_{j=1}^{N_v} 
  \begin{pmatrix}[1.2]  1 \\ v_j \\ v_j^2 \end{pmatrix} {\mathcal M}^n_{ij} + \left( 1 - \theta \right) \sum_{j=1}^{N_v} \begin{pmatrix}[1.2] 1  \\ v_j \\ v_j^2 \end{pmatrix} f^n_{ij}.
\end{equation}
In general, mass, momentum, and energy will not be conserved:
\begin{equation}
\sum_{j=1}^{N_v} \begin{pmatrix}[1.2] 1  \\ v_j \\ v_j^2 \end{pmatrix} f^{n+1}_{ij} \, \, \ne \, \, \sum_{j=1}^{N_v} \begin{pmatrix}[1.2] 1  \\ v_j \\ v_j^2 \end{pmatrix} f^n_{ij},
\end{equation}
because the distribution, $f$, and the Maxwell-Boltzmann distribution, ${\mathcal M}$, do not have exactly matching discrete moments:
\begin{equation}
\label{eqn:non-matching-moments}
\sum_{j=1}^{N_v} \begin{pmatrix}[1.2] 1  \\ v_j \\ v_j^2 \end{pmatrix} f^{n}_{ij} \, \, \ne \, \, \sum_{j=1}^{N_v} \begin{pmatrix}[1.2] 1  \\ v_j \\ v_j^2 \end{pmatrix} {\mathcal M}^n_{ij}.
\end{equation}
There are two reasons that the moments in \eqref{eqn:non-matching-moments} do not match:
\begin{enumerate}
\item The velocity extent is limited to the domain $v \in \left[v_{\text{min}}, v_{\text{max}}\right]$, but the Maxwell-Boltzmann distribution satisfies ${\mathcal M} > 0$ for all $v \in \left(-\infty, \infty\right)$, which means that any portion of the distribution that contributes to the total mass, momentum, or energy for $v>v_{\text{max}}$ or $v<v_{\text{min}}$ is lost to the numerical method; and
\item The numerical quadrature is only approximate and cannot exactly integrate the Maxwell-Boltzmann distribution on a finite mesh.
\end{enumerate}

We propose to overcome limitation \eqref{eqn:non-matching-moments} by replacing the discrete Maxwell-Boltzmann distribution, ${\mathcal M}_{ij}$, by a modified version, $\widetilde{\mathcal M}_{ij}$, that will allow us to recover
exact conservation. The idea is similar in spirit to Mieussens \cite{Mieussens2000DiscreteDynamics}, but we propose a more straightforward approach that only requires the solution of a linear system, and the solution can be written down in closed form.
We consider a modified Maxwell-Boltzmann distribution of the form:
\begin{equation}
\label{eqn:mod_maxwellian}
	\widetilde{\mathcal M}_{ij} = \frac{\rho_i}{\sqrt{2\pi T_i}} \exp\left[-\frac{\mu_{ij}^2}{2} \right]\Bigl( a_{1i} + a_{2i} \, \mu_{ij} + a_{3i}  \left( \mu_{ij}^2-1 \right) \Bigr) \quad \text{where} \quad
	\mu_{ij} = \frac{v_j-u_i}{\sqrt{T_i}}.
\end{equation}
This form also borrows from the classical Grad moment-closure \cite{Grad1949OnGases} in that 
\begin{equation}
\left\{ 1, \, \mu, \, \mu^2-1 \right\} 
\end{equation}
are the first three Hermite polynomial basis functions (i.e., orthogonal polynomials on $\mu\in\left(-\infty,\infty\right)$ with weight $\exp\left(-\mu^2/2 \right)$). If there were no discrete errors in the numerical integration, the coefficients would be $a_{1i}=1$, $a_{2i}=0$, and $a_{3i}=0$, which would recover the original Maxwell-Boltzmann distribution, ${\mathcal M}_{ij}$. 

\begin{algorithm}[!t]
\caption{Initialize computation \label{alg:alg-init}}
\begin{algorithmic}
\STATE {\bf INPUT:} \hspace{6.5mm} $N_x$, $N_v$, $x_{\text{low}}$, $x_{\text{high}}$,
	$v_{\text{low}}$, $v_{\text{high}}$, $f_0(x,v)$ (PDF initial condition function)
\STATE {\bf OUTPUT:} \hspace{2mm} $f_{ij}$ (PDF at initial time), $x_i$, $v_j$, $\Delta x$, $\Delta v$, $v_{\text{max}}$
\STATE
\STATE $\Delta x = \left( x_{\text{high}} - x_{\text{low}} \right)/N_x$ \hspace{10mm} \inmag{\tt // Grid spacing}
\STATE $\Delta v = \left( v_{\text{high}} - v_{\text{low}} \right)/N_v$ \hspace{11mm} \inmag{\tt // Velocity grid spacing}
\STATE $v_{\text{max}} = \text{max}\left\{ \bigl|v_{\text{low}}\bigr|, \bigl|v_{\text{high}}\bigr| \right\}$ \hspace{4mm} \inmag{\tt // Maximum particle velocity}
\STATE
\FOR{$i = 1,\ldots,N_x$}
\STATE $x_i = x_{\text{low}} + \left( i - 0.5 \right) \Delta x$ \hspace{6mm} \inmag{\tt // x cell center}
\ENDFOR
\STATE
\FOR{$j = 1,\ldots,N_v$}
\STATE $v_j = v_{\text{low}} + \left( j - 0.5 \right) \Delta v$ \hspace{6mm} \inmag{\tt // v cell center}
\ENDFOR
\STATE
\FOR{$i = 1,\ldots,N_x$}
\FOR{$j = 1,\ldots,N_v$}
\STATE $f_{ij} = f_0\left(x_i, v_j \right)$ \hspace{17.5mm} \inmag{\tt // initial PDF at cell centers}
\ENDFOR
\ENDFOR
\end{algorithmic}
\end{algorithm}

To compute the unknown coefficients, $a_{1i}$, $a_{2i}$, and $a_{3i}$, we take the modified collision sub-step update:
\begin{equation}
  f_{ij} \leftarrow \theta \widetilde{\mathcal M}_{ij} + \left( 1 - \theta \right) f_{ij},
\end{equation}
and require that $\widetilde{\mathcal M}$ and $f$ have the same {\it discrete} mass, momentum, and energy moments:
\begin{equation}
\frac{\Delta v}{\sqrt{2\pi T_i}} \sum_{j=1}^{N_v} \begin{pmatrix}[1.2] 1  \\ \mu_{ij} \\ \mu_{ij}^2 \end{pmatrix} \exp\left[-\frac{\mu_{ij}^2}{2} \right]\Bigl( a_{1i} + a_{2i} \, \mu_{ij} + a_{3i}  \left( \mu_{ij}^2-1 \right) \Bigr) = 
\frac{\Delta v}{\rho_i} \sum_{j=1}^{N_v} \begin{pmatrix}[1.2] 1  \\ \mu_{ij} \\ \mu_{ij}^2 \end{pmatrix} f_{ij},
\end{equation}
which we can write as a $3\times3$ linear system:
\begin{equation}
\begin{bmatrix}[1.2]
A_{0i} & A_{1i} & A_{2i}-A_{0i} \\
A_{1i} & A_{2i} & A_{3i}-A_{1i} \\
A_{2i} & A_{3i} & A_{4i}-A_{2i}
\end{bmatrix}
\begin{bmatrix}[1.2]
  a_{1i} \\ a_{2i} \\ a_{3i}
\end{bmatrix} =
\begin{bmatrix}[1.2]
  1 \\ 0 \\ 1
\end{bmatrix},
\end{equation}
where
\begin{equation}
\label{eqn:cons_coeff_1}
A_{ki} = \frac{\Delta v}{\sqrt{2 \pi T_i}} \sum_{j=1}^{N_v} \mu_{ij}^k \, \exp\left[-\frac{\mu_{ij}^2}{2} \right]
\quad \text{for} \quad k=0,1,2,3,4.
\end{equation}
This $3\times3$ linear system can be easily solved exactly:
\begin{gather}
\label{eqn:cons_coeff_2}
	\begin{bmatrix}
	a_{1i} \\ a_{2i} \\ a_{3i}
	\end{bmatrix} = 
	\frac{1}{d_i} 
	\begin{bmatrix}
       A_{1i}^2 + A_{2i} \left(2 A_{2i} - A_{0i} - A_{4i} \right) - A_{3i} \left(2 A_{1i} - A_{3i} \right) \\
       A_{1i} \left( A_{4i} - A_{2i} \right) + A_{3i} \left( A_{0i} - A_{2i} \right) \\
       A_{1i} \left( A_{1i} - A_{3i} \right) + A_{2i} \left( A_{2i} - A_{0i} \right)
    \end{bmatrix},
\end{gather}
where
\begin{gather}
\label{eqn:cons_coeff_3}
	d_i = A_{2i}^3 - 2 A_{1i} A_{2i} A_{3i} + A_{0i} A_{3i}^2 + A_{1i}^2 A_{4i} - A_{0i} A_{2i} A_{4i}.
\end{gather}
We now have a suitable modified Maxwell-Boltzmann distribution, $\widetilde{\mathcal M}_{ij}$, as defined in \eqref{eqn:mod_maxwellian}. 

To finalize the description of the proposed numerical method for the 1D1V Boltzmann-BGK equation \eqref{eqn:Boltzmann-BGK}, we summarize the scheme in two short chunks of pseudocode. In Algorithm \ref{alg:alg-init}, we show how to initialize the mesh and the distribution function. In Algorithm \ref{alg:alg-timestep}, we summarize the complete time-stepping procedure, including how the operator splitting can be implemented. 
We implement the above-described algorithm in a freely available code \cite{code:VBRossmanith2025}; the central portion of the code is written in Java, and the plotting routines are written in Python using NumPy and Matplotlib.

\begin{algorithm}[!t]
\caption{Advance PDF and moments from initial condition to final time \label{alg:alg-timestep}}
\begin{algorithmic}
\STATE {\bf INPUT:} \hspace{6.5mm} $f_{ij}$ (initial PDF), $v_j$, $\Delta x$, $\Delta v$, $N_x$, $N_v$, $\text{CFL}$, $v_{\text{max}}$, $\varepsilon$, $\text{finalTime}$
\STATE {\bf OUTPUT:} \hspace{2mm} $f_{ij}$ (PDF at final time)
\STATE
\STATE $\Delta t = \Delta x \, {\text{CFL}}\bigl/{v_{\text{max}}}$ \hspace{12mm} \inmag{\tt // Set time step}
\STATE $\text{nSteps} = \lceil\text{finalTime}\big/{\Delta t}\rceil$ \hspace{3mm}  \inmag{\tt // Compute number of time steps}
\STATE  $\Delta t = \text{finalTime}\bigl/{\text{nSteps}}$ \hspace{6mm} \inmag{\tt // Reset time step}
\STATE  $\text{CFL} = v_{\text{max}} \, {\Delta t}\bigl/{\Delta x}$ \hspace{12mm} \inmag{\tt // Reset CFL number}
\STATE  $\theta = \frac{\Delta t \left(\Delta t + 48 \varepsilon\right)}{\left(\Delta t + 6 \varepsilon\right) \left(\Delta t + 8 \varepsilon\right)}$ \hspace{15mm} \inmag{\tt // Set collision parameter}
\STATE
\FOR{$n = 1,\ldots,\text{nSteps}$}
    \STATE {\bf Step 1:} From $f_{ij}$ compute moments via \eqref{eqn:discrete_moments} and $a_{1i}$,
    $a_{2i}$, $a_{3i}$ via \eqref{eqn:cons_coeff_1}--\eqref{eqn:cons_coeff_3}
    \STATE {\bf Step 2:} $\Delta t/2$ collision step: $f_{ij} \leftarrow \theta \widetilde{\mathcal M}_{ij} + \left(1-\theta\right) f_{ij}$ using \eqref{eqn:mod_maxwellian}
    \STATE {\bf Step 3:} Copy $f_{ij}$ into a temporary array:  $f^{\star}_{ij} \leftarrow f_{ij}$
    \STATE {\bf Step 4:} $\Delta t$ transport step: compute new $f_{ij}$ from old $f^{\star}_{ij}$ using \eqref{eqn:transport1} and \eqref{eqn:transport2}
    \STATE {\bf Step 5:} From $f_{ij}$ compute moments via \eqref{eqn:discrete_moments} and $a_{1i}$,
    $a_{2i}$, $a_{3i}$ via \eqref{eqn:cons_coeff_1}--\eqref{eqn:cons_coeff_3}
    \STATE {\bf Step 6:} $\Delta t/2$ collision step: $f_{ij} \leftarrow \theta \widetilde{\mathcal M}_{ij} + \left(1-\theta\right) f_{ij}$ using \eqref{eqn:mod_maxwellian}
\ENDFOR
\end{algorithmic}
\end{algorithm}


\begin{figure}[!t]
\begin{tabular}{cc}
(a)\includegraphics[width=8cm]{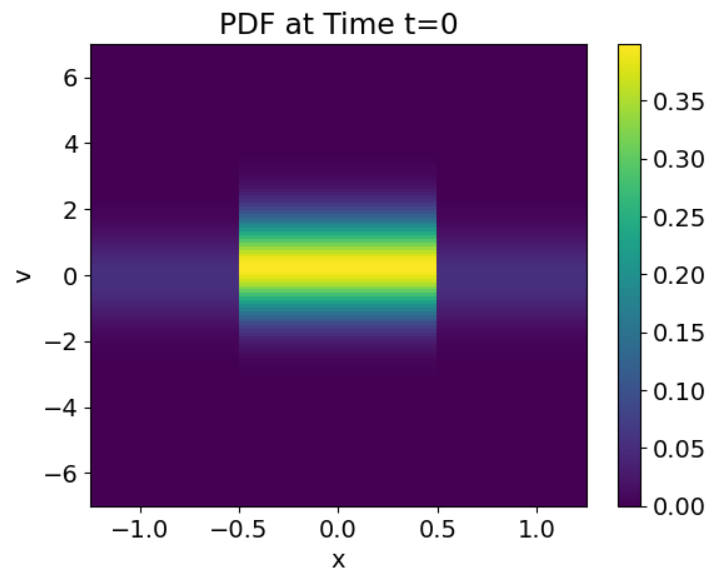} &
(b)\includegraphics[width=8cm]{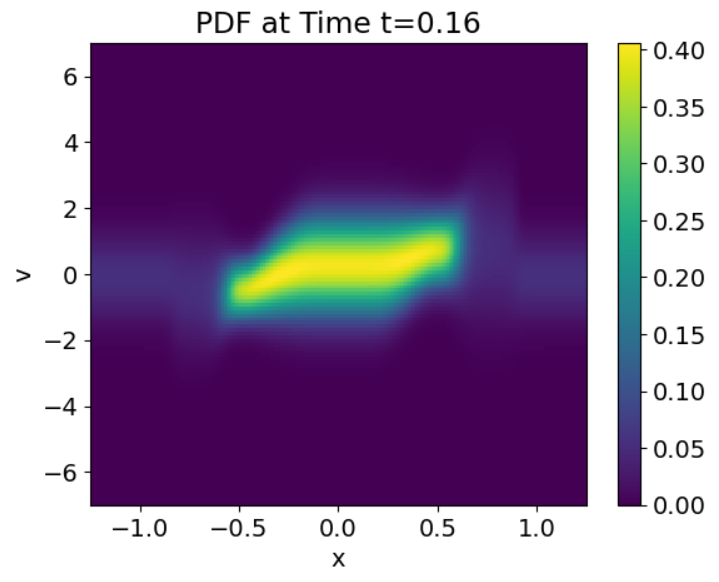}
\end{tabular}
\caption{Particle distribution function in the $(x,v)$ phase space at times (a) $t=0$ and (b) $t=0.16$. We used the mesh and parameter values described in \S\ref{sec:numerical-results}.\label{fig:pdf}}
\end{figure}

\section{Numerical Results}
\label{sec:numerical-results}
To demonstrate the effectiveness of the modified collision operator, we consider a test case with the following initial data:
\begin{gather}
f\left(t=0,x,v\right) = {\mathcal M}_0\left(x,v\right) = \frac{\rho_0(x)}{\sqrt{2\pi  T_0(x)}} \exp\left[-\frac{\left(v-u_0(x)\right)^2}{2  T_0(x)} \right],
\end{gather}
where
\begin{gather}
\left( \rho_0, \, u_0, \, T_0 \right)\left(x\right) = 
   \begin{cases}
      \left( 1.000, \, 0.250, \, 1.000 \right) & \quad \text{if} \, \left| x \right| < 0.5, \\
      \left( 0.125, \, -0.10, \, 0.800 \right) & \quad \text{otherwise},
   \end{cases}
\end{gather}
on the domain $x\in\left[-1.25, 1.25 \right]$ with {\it periodic} boundary conditions:
\begin{equation}
f\left(t,x=-1.25,v\right) = f\left(t,x=+1.25,v\right).
\end{equation}
We choose the truncated velocity domain $v \in [-7,7]$ and the following parameters:
\begin{equation}
N_x = 256, \quad N_v = 128, \quad \varepsilon = 0.01, \quad \text{CFL} = 1.95, \quad T_{\text{final}}  = 0.16.
\end{equation}
The periodic boundary conditions are enforced during the transport sub-step via \eqref{eqn:transport1} and \eqref{eqn:transport2}
by replacing all indeces that stray out of the range $i=1,2,\ldots,N_x$ by the following:
\begin{equation}
\left(-1\right) \rightarrow \left(N_x-1\right), \quad
\left(0\right) \rightarrow \left(N_x\right), \quad
\left(N_x+1\right)\rightarrow \left(1\right), \quad
\left(N_x+2\right)\rightarrow \left(2\right).
\end{equation}

\begin{figure}[!t]
\begin{tabular}{cc}
(a)\includegraphics[width=7.5cm]{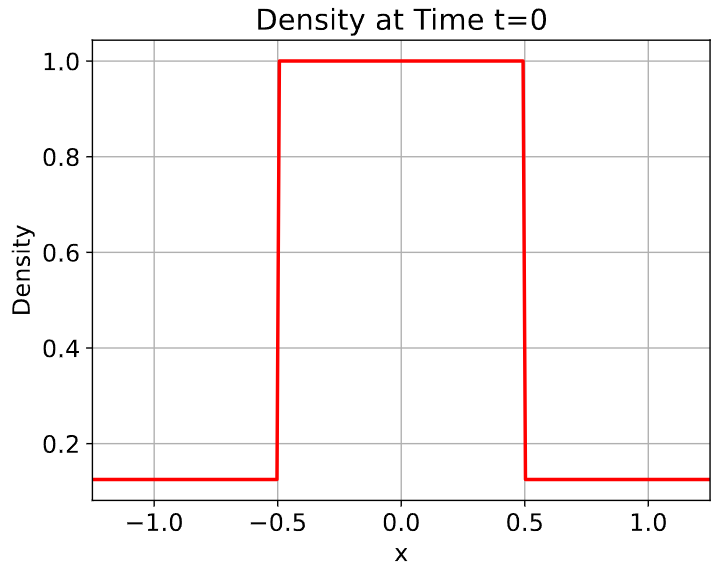} &
(b)\includegraphics[width=7.5cm]{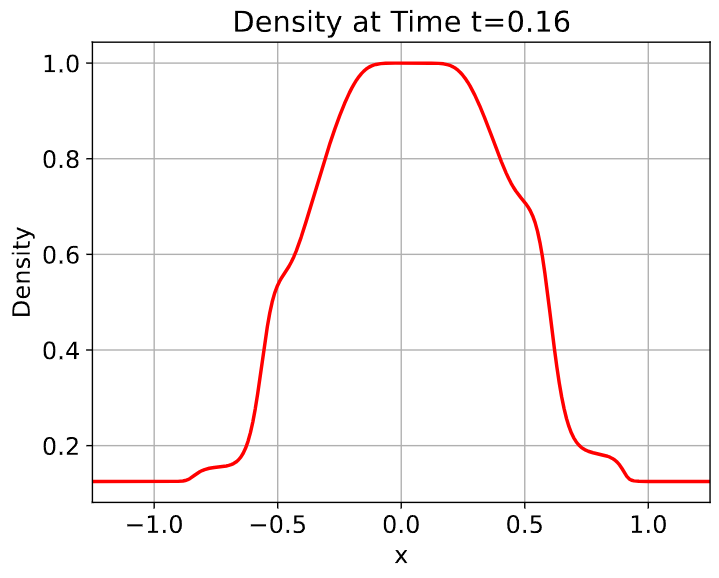} \\
(c)\includegraphics[width=7.5cm]{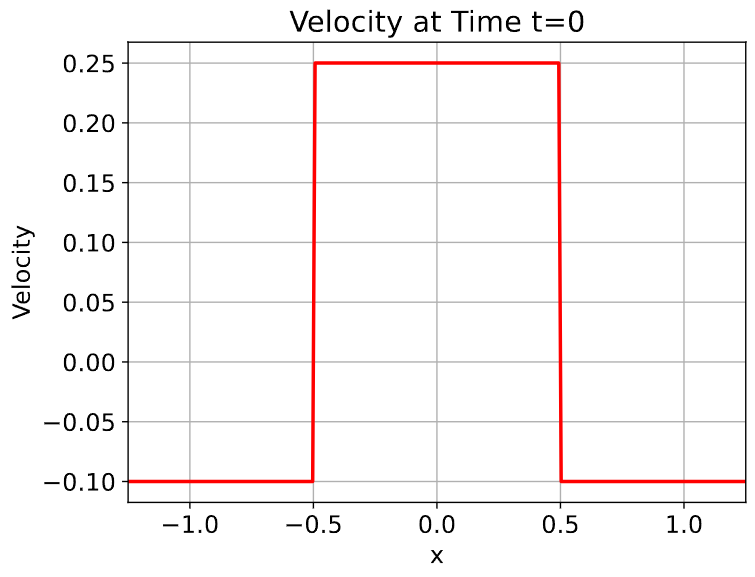} &
(d)\includegraphics[width=7.5cm]{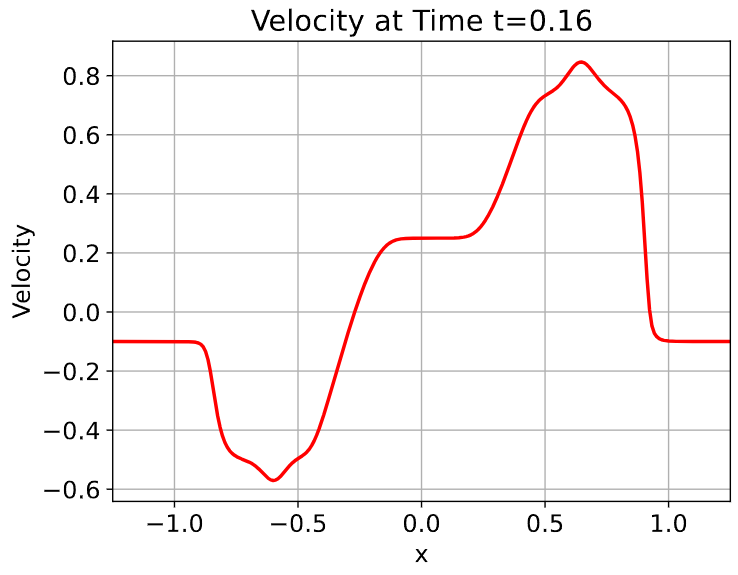} \\
(e)\includegraphics[width=7.5cm]{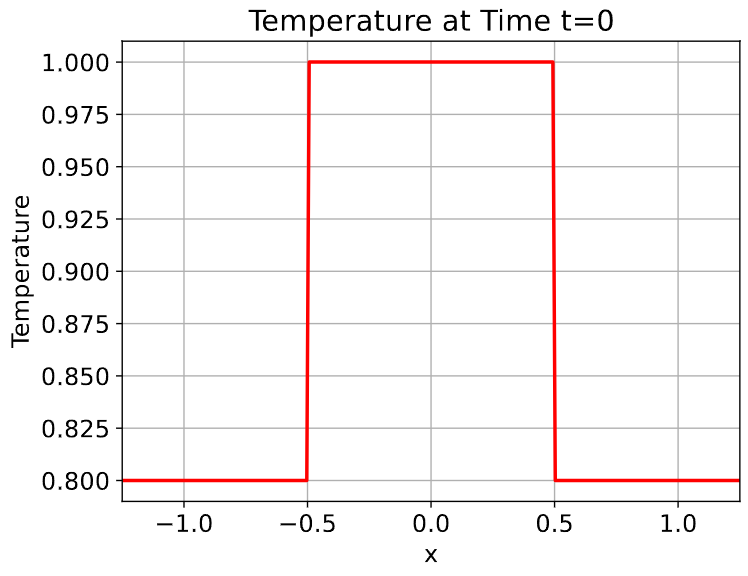} &
(f)\includegraphics[width=7.5cm]{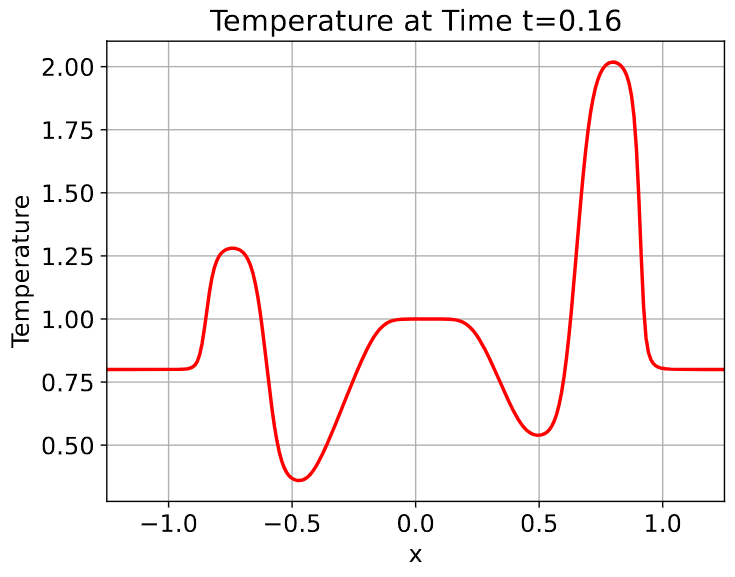}
\end{tabular}
\caption{Plots of various moments of the distribution function. Shown in these panels are the 
(a) density ($\rho$) at time $t=0$, (b) density ($\rho$) at time $t=0.16$,
(c) fluid velocity ($u$) at time $t=0$, (d) fluid velocity ($\rho$) at time $t=0.16$,
(e) temperature ($T$) at time $t=0$, and (f) temperature ($T$) at time $t=0.16$. We used the mesh and parameter values described in \S\ref{sec:numerical-results} \label{fig:moments}.}
\end{figure}

We implemented the algorithms described in the previous section, Algorithms \ref{alg:alg-init} and \ref{alg:alg-timestep}, in freely available Java code, with corresponding plotting routines in Python (via NumPy and Matplotlib). The particle density function (PDF) in the $(x,v)$ phase space is shown at the initial ($t=0$) and final ($t=0.16$) times in Figure \ref{fig:pdf}. Corresponding plots of the various moments of the distribution function are shown in Figure \ref{fig:moments}; these panels display the (a) density ($\rho$) at time $t=0$, (b) density ($\rho$) at time $t=0.16$, (c) fluid velocity ($u$) at time $t=0$, (d) fluid velocity ($\rho$) at time $t=0.16$, (e) temperature ($T$) at time $t=0$, and (f) temperature ($T$) at time $t=0.16$. These plots show that the initial discontinuous block in the center of the domain breaks up into shocks and rarefactions that propagate to the left and right. Because the Knudsen number is set to $\varepsilon=0.01$, we also see the effect of viscosity on smoothing out shocks and rarefactions. These simulation results are consistent with similar problems reported in the literature, including in Sar \cite{thesis:Sar2024} and Pieraccini and Puppo \cite{Pieraccini2007Implicit-explicitEquations}. 

\begin{figure}[!t]
\begin{tabular}{cc}
(a)\includegraphics[width=8cm]{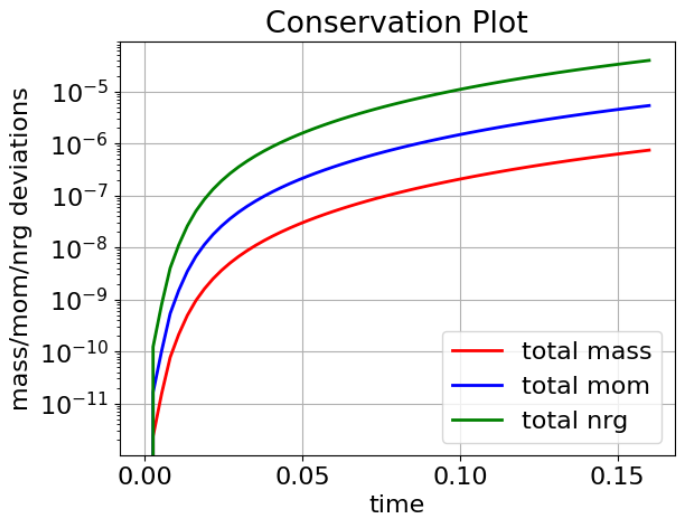} &
(b)\includegraphics[width=8cm]{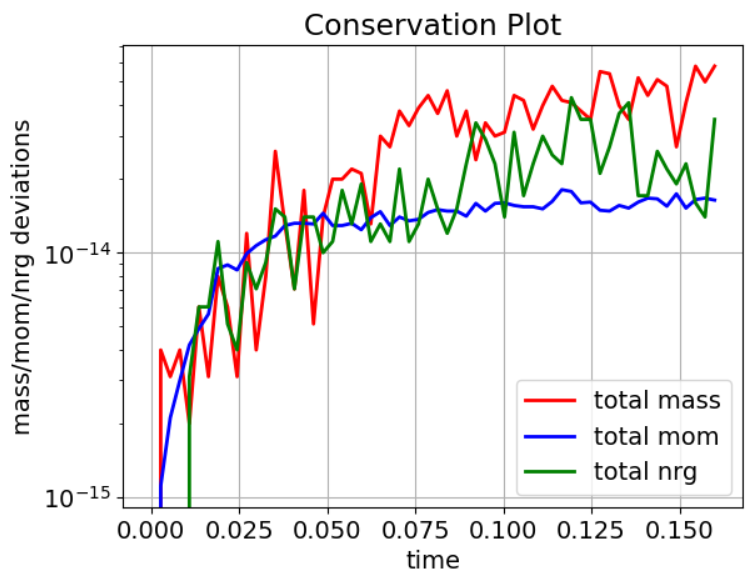}
\end{tabular}
\caption{Conservation of mass, momentum, and energy. Shown in Panel (a) is the simulation without the collision operator modification (i.e., ${\mathcal{M}}\leftarrow\widetilde{\mathcal{M}}$). Shown in Panel (b) is the simulation with the proposed modified collision operator, $\widetilde{\mathcal{M}}$. These plots clearly show the lack of exact conservation when using ${\mathcal{M}}$; furthermore, the conservation errors tend to increase over time. The plots also clearly show the recovery of exact conservation (at least up to machine precision) with the proposed modified collision operator. \label{fig:conservation}}
\end{figure}

The results shown here are with the modified Maxwell-Boltzmann distribution, $\widetilde{\mathcal{M}}$, as defined by \eqref{eqn:mod_maxwellian}. We have also run this example on the identical method, but with the standard Maxwell-Boltzmann distribution, $\mathcal{M}$, and, at least on the scale of the plots, the results with either $\widetilde{\mathcal{M}}$ or ${\mathcal{M}}$ are essentially indistinguishable. Upon closer inspection, however, there is a key difference: the modified approach with $\widetilde{\mathcal{M}}$ exactly conserves -- up to machine precision -- the mass, momentum, and energy, while the same method with ${\mathcal{M}}$ does not. We quantify this by defining the deviations from the initial total mass, momentum, and energy:
\begin{align}
  \text{relative total mass change:} & \quad \Delta \rho^n := \left| \sum_{i=1}^{N_x} \sum_{j=1}^{N_v} \left( f^{n}_{ij} - f^{0}_{ij}\right) \right|\Biggl/ \left| \sum_{i=1}^{N_x} \sum_{j=1}^{N_v} f^{0}_{ij} \right|, \\
  \text{relative total momentum change:} & \quad \Delta m^n := \left| \sum_{i=1}^{N_x} \sum_{j=1}^{N_v} \left( v_j f^{n}_{ij} - v_j f^{0}_{ij}\right) \right|\Biggl/ \left| \sum_{i=1}^{N_x} \sum_{j=1}^{N_v} v_j f^{0}_{ij} \right|, \\
  \text{relative total energy change:} & \quad \Delta {\mathcal E}^n := \left| \sum_{i=1}^{N_x} \sum_{j=1}^{N_v} \left( v_j^2 f^{n}_{ij} - v_j^2 f^{0}_{ij}\right) \right|\Biggl/ \left| \sum_{i=1}^{N_x} \sum_{j=1}^{N_v} v^2_j f^{0}_{ij} \right|.
\end{align}
The relative deviations are shown in Figure \ref{fig:conservation}, where the lack of exact conservation when using ${\mathcal{M}}$ is shown in \ref{fig:conservation}(a), while exact preservation (at least up to machine precision) with the proposed modified collision operator is shown in \ref{fig:conservation}(b). We also note that the relative deviations accumulate and grow over time when using ${\mathcal{M}}$.

Finally, we note that the Java/Python code that produced the plots in this section can be freely downloaded and executed \cite{code:VBRossmanith2025}; the central portion of the code is written in Java, and the plotting routines are written in Python using NumPy and Matplotlib.


\section{Conclusion}
\label{sec:conclusion}
In this work, we developed a novel numerical method for solving the 1D1V Boltzmann-BGK equation. Several important ingredients were combined to create the resulting numerical method. First, we made use of Strang operator splitting \cite{article:St68} to create separate transport and collision sub-steps, each of which was easier to discretize than the whole system; second, we introduced a third-order accurate Lax-Wendroff-type scheme \cite{Lax1960SystemsLaws} for the transport sub-step; third, we made use of the second-order L-stable TR-BDF method \cite{article:TRBDF2} for the collision sub-step. These ingredients resulted in an accurate, efficient, and robust numerical method across all Knudsen numbers, $\varepsilon>0$, including in the singular limit: $\varepsilon \rightarrow 0^+$. The final component we introduced is a novel treatment of the collision sub-step to guarantee that the total mass, momentum, and energy are conserved even in the presence of a truncated velocity range and quadrature errors in computing moments. The key idea was to multiply the Maxwell-Boltzmann distribution in the collision sub-step by a quadratic Hermite polynomial, where the coefficients are chosen to ensure exact conservation. The resulting scheme was verified on piecewise constant initial data with periodic boundary conditions; exact conservation up to machine precision was demonstrated. The resulting Java/Python code is freely available \cite{code:VBRossmanith2025}; the central portion of the code is written in Java, and the plotting routines are written in Python using NumPy and Matplotlib. Future extensions of this work will focus on the higher-dimensional variants of the Boltzmann-BGK equation.

\bibliographystyle{siamplain}


\end{document}